\magnification=1200

\advance\baselineskip by 0.5pt
\font\rm=pplr at 10pt
\font\it=pplri at 10pt
\font\bf=pplb at 10pt

\font\sc=pplr at 7pt
\rm

\footline={\hss\rm\folio\hss}
\parskip=10pt plus3pt minus3pt
\def\itc#1{{\it #1\/}}
\newcount\secno \secno=0
\def\sec#1\par{\global\advance\secno by 1\bigskip
    \noindent{\bf\the\secno. #1}\par\nobreak\noindent}

\centerline{\bf The Surprise Examination or Unexpected Hanging Paradox}
\medskip
\centerline{Timothy Y. Chow}
\medskip
\centerline{[As published in \itc{Amer.\ Math.\ Monthly}
{\bf 105} (1998), 41--51; ArXived with permission.]}
\medskip

Many mathematicians have a dismissive attitude towards paradoxes.
This is unfortunate, because many paradoxes are rich in content,
having connections with serious mathematical ideas as well as
having pedagogical value in teaching elementary logical reasoning.
An excellent example is the so-called ``surprise examination
paradox'' (described below),
which is an argument that seems at first to be too silly
to deserve much attention.  However, it has inspired an
amazing variety of philosophical and mathematical investigations that
have in turn uncovered links to G\"odel's incompleteness theorems,
game theory, and several other logical paradoxes (e.g., the liar
paradox and the sorites paradox).  Unfortunately, most mathematicians
are unaware of this because most of the literature has been published
in philosophy journals.

In this article, I describe some of this work,
emphasizing the ideas that are particularly interesting mathematically.
I also try to dispel some of the confusion
that surrounds the paradox and
plagues even the published literature.
However, I do not try to correct every error or explain every idea
that has ever appeared in print.
Readers who want more comprehensive surveys should see
[30, chapters 7 and~8], [20], and [16].

At times I assume some knowledge of mathematical logic
(such as may be found in Enderton~[10]), but the reader who lacks
this background may safely skim these sections.

\sec The paradox and the meta-paradox

Let us begin by recalling the paradox.   It has many variants,
the earliest probably being Lennart Ekbom's surprise drill,
and the best known to mathematicians (thanks to Quine and Gardner)
being an unexpected hanging.
We shall give the surprise examination version.

{\it\narrower\noindent
A teacher announces in class that an examination will be held on some
day during the following week, and moreover that the examination will
be a surprise.  The students argue that a surprise exam cannot occur.
For suppose the exam were on the last day of the week.  Then on the
previous night, the students would be able to predict that the exam
would occur on the following day, and the exam would not be a surprise.
So it is impossible for a surprise exam to occur on the last day.  But
then a surprise exam cannot occur on the penultimate day, either, for
in that case the students, knowing that the last day is an impossible
day for a surprise exam, would be able to predict on the night before
the exam that the exam would occur on the following day.  Similarly,
the students argue that a surprise exam cannot occur on any other day
of the week either.  Confident in this conclusion, they are of course
totally surprised when the exam occurs (on Wednesday, say).  The
announcement is vindicated after all.  Where did the students'
reasoning go~wrong?\par}

The natural reaction to a paradox like this is
to try to resolve it.  Indeed, if you have not seen this paradox
before, I encourage you to try to resolve it now before reading on.
However, I do not want to
discuss the resolution of the paradox right away.
Instead, for reasons that should become apparent,
I discuss what I call the ``meta-paradox'' first.

The meta-paradox consists of two seemingly incompatible facts.  The
first is that the surprise exam paradox seems easy to resolve.
Those seeing it for the first time typically have the instinctive
reaction that the flaw in the students' reasoning is obvious.
Furthermore, most readers who have tried to think it through
have had little difficulty resolving it to their own satisfaction.

The second (astonishing) fact is that
to date nearly a hundred papers on the paradox have been published, and
still no consensus on its correct resolution has been reached.  The
paradox has even been called a ``significant problem'' for philosophy
[30, chapter~7, section~VII].  How can this be?  Can such a
ridiculous argument really be a major unsolved mystery?  If not,
why does paper after paper begin by brusquely dismissing all previous work
and claiming that it alone presents the long-awaited simple solution
that lays the paradox to rest once and for all?

Some other paradoxes suffer from a similar meta-paradox,
but the problem is especially acute in the case of the surprise
examination paradox.  For most other trivial-sounding paradoxes
there is broad consensus on the proper resolution,
whereas for the surprise exam paradox there is not even
agreement on its proper formulation.
Since one's view of the meta-paradox influences the way one views
the paradox itself, I must try to clear up the former before
discussing the latter.

In my view, most of the confusion has been caused by authors who
have plunged
into the process of ``resolving'' the paradox without first having
a clear idea of what it \itc{means} to ``resolve'' a paradox.
The goal is poorly understood, so controversy over whether the
goal has been attained is inevitable.
Let me now suggest a way of thinking about the process
of ``resolving a paradox'' that I believe dispels the meta-paradox.

In general, there are two steps involved in resolving a paradox.
First, one establishes precisely \itc{what the paradoxical argument~is}.
Any unclear terms are defined carefully
and all assumptions and logical steps are stated clearly and explicitly,
possibly in a formal language of some kind.
Second, one \itc{finds the fault in the argument.}
Sometimes, simply performing step one reveals the flaw, e.g., when
the paradox hinges on confusing two different meanings of the same word,
so that pointing out the ambiguity suffices to dispel the confusion.
In other cases, however,
something more needs to be done; one must locate the bad assumptions,
the bad reasoning, or (in desperate circumstances)
the flaw in the structure of logic itself.

These two steps seem straightforward, but there are a few subtleties.
For example, if, in the second step, the flaw is caused by bad assumptions,
it may be hard to isolate a unique culprit.
Sometimes what we discover is a set of mutually incompatible assumptions
such that rejecting any one of them suffices to eliminate the contradiction.
When this occurs, however, notice that
while it may be an interesting question to decide which assumption to reject,
such a decision is \itc{not} usually needed to resolve the paradox.
It is usually enough to exhibit the incompatible assumptions
and state that their joint inconsistency is the source of the paradox.

The first step of resolving a paradox can also be subtle.
As many investigators of the surprise exam paradox have noted,
formal versions of a paradox sometimes miss the essence of
the original informal version.
Such a mistranslation evades the paradox instead of resolving it.
Certainly, this is a real danger,
and numerous authors have fallen into this trap.
However,
there is a simple but important point here that is often overlooked:
the question of whether or not a particular
formalization of a paradox ``captures its essence''
is to some extent a matter of opinion.
Given two formalizations of the paradox,
one person may think that the first captures the essence better
but another may prefer the second.
One cannot say who is objectively right, since
there is always some vagueness in the original informal account.
To be sure, one can sometimes argue that a particular formalization is
inadequate by proposing a variation of the paradox that seems to retain
its essence but for which the particular formalization fails.
Even here, though, there is some room for differences of opinion,
because one can sometimes argue that the variant paradox does not in fact
retain the essence but is actually a different paradox that requires
a different solution.

Thus, sometimes there exist multiple formalizations
of a paradox that all capture its essence reasonably well.
In such cases I believe it is misguided to speak of \itc{the} resolution
of the paradox.  This point has also been made by Kirkham~[16].

With these ideas in mind we can easily explain the meta-paradox.
A careful look at the literature confirms our suspicion
that the paradox is not hard to resolve,
because most authors have succeeded in finding resolutions.
Most of the controversies have been false controversies.
For example, there has been much debate between what I
call the ``epistemological school'' (which formalizes the paradox
using concepts such as knowledge, belief and memory) and the
``logical school'' (which avoids such concepts) over who has
the ``right'' formalization.  But
both approaches are reasonable
and neither is guilty of evasion.

Also, within the epistemological school there has been
much debate over which axiom of a certain set of mutually inconsistent
axioms about knowledge should be rejected.  The question is an
interesting one from the point of view of philosophically analyzing the
concept of knowledge, but if we agree that identifying the
``right'' axiom to reject is not essential to resolving the paradox
then this debate need not trouble us.

Having dealt with the meta-paradox, we now turn to the paradox
itself and explore several different approaches.

\sec The logical school

We mathematicians have a firm belief that logic and mathematics are
consistent.  When we are confronted with a paradox, therefore, our
tendency is to assume, even before analyzing the paradox, that either
the paradox cannot be translated into a purely logical or mathematical
argument, or that if it can be so translated, the faulty step or
assumption will become immediately apparent.  So a natural reaction
to the surprise examination paradox (at least for a mathematician)
is to take the students' argument and try to convert it into a
rigorous proof in order to find the flaw.  Let us now do
this and see what happens.

Every proof begins with axioms.  The students' argument seems to deduce
a contradiction from the teacher's announcement, so it seems that the
axioms in this case ought to be some formalization of the announcement.
Now, part of the announcement---the claim that an examination will take
place some time during the following week---is not difficult to formalize,
but the part that says that the examination will be a surprise is not as
clear.  What is meant by ``surprise''?

Whatever ``surprise'' means, it must at least mean that the students will
not be able to deduce logically the date of the examination ahead of time,
for if the students could \itc{prove} that the date of the examination
were such-and-such before the date arrived, they would surely not be at
all surprised by the exam.  So a first step towards formalizing the
teacher's announcement might be, ``There will be an examination next
week and its date will not be deducible in advance.''

This is not sufficient, however, because every proof
begins with axioms.  To say that the date of the examination will not be
deducible in advance is a vague statement until the axioms from which the
date cannot be deduced are specified precisely.  Now, it is not completely
clear which axioms are in question here; the informal word ``surprise'' is
too vague to give us many clues.  However, if our formalization is to be
at all true to the original paradox, it should at least allow us to
formalize the students' argument to some degree.
Formalizing the teacher's announcement as

\item{(a)}
There will be an examination
next week and its date will not be deducible in advance
from an empty set of assumptions

\noindent
is certainly not satisfactory because it does not allow the students'
argument even to begin.  This would evade the paradox and not resolve it.

A better attempt at formalization might be something like

\item{(b)} There will be
an examination next week and its date will not be deducible
in advance from the assumption that the examination will occur
some time during the week.

This formalization allows at least the first
step of the students' argument to be carried out: given this announcement,
the students can deduce that the examination will not occur on the
last day of the week.  However, if we try to reproduce the next step of the
students' argument---the step that eliminates the penultimate day of the
week---we find ourselves stuck.  In order to eliminate the penultimate day,
the students need to argue that their ability to deduce, from statement~(b),
that the examination will not occur on the last day implies that a last-day
examination \itc{will not be surprising.}  But since we have restricted
``surprising'' to mean ``not deducible from the assumption that the
examination will occur sometime during the week'' instead of ``not deducible
from \itc{statement~(b),}'' the students' argument is blocked.
To continue the argument we need to be able to use the nondeducibility
from the announcement as an assumption, i.e., we must embed the
nondeducibility from the announcement into the announcement itself.

It now becomes clear that to carry out the students' argument,
one needs a formalization that is something like

\item{(c)} There will be an examination next week and its date will
not be deducible in advance using \itc{this announcement} as an axiom.

In other words, the announcement must be formulated as a
\itc{self-referential} statement!

There is a temptation to end the analysis here with a comment
that the self-referential nature of statement~(c)
is the source of the paradox.
After all, if a self-referential definition like this were to be
presented in a mathematical paper, we would surely reject it
instantly as illegal.
Indeed, Shaw concludes his paper~[25] with just such a comment.

However, we need to be careful.
For it \itc{is} possible in mathematics to formalize
certain kinds of self-referential statements.
Indeed, this was one of the crucial ideas
in G\"odel's proof of his incompleteness theorems,
and it is now a standard technique in mathematical logic.
It is natural to ask if this technique can be used to obtain a
completely formal version of statement~(c).
The answer is yes; we give the construction
(due to Fitch~[11]; [4] and~[32] have similar constructions)
in some detail since it is rather interesting.

\def\Do{\itc{Q$_1$}}
\def\Dt{\itc{Q$_2$}}
Let us reduce the number of days to two for simplicity
(we consider one-day weeks shortly),
and let \Do\ and \Dt\ be statements representing the occurrence
of the exam on days one and two, respectively.
Then what we are seeking is a statement~\itc{S} such that

\bigskip
\centerline{\itc{S} $\equiv$ (\Do\ \& ([\itc{S} $\Rightarrow$ \Do] is
unprovable)) or else (\Dt\ \& ([\itc{S} \& $\sim$\Do\ $\Rightarrow$
\Dt] is unprovable)).}

\noindent
Given a first-order language that contains enough elementary
arithmetic to handle primitive recursive functions,
together with some G\"odel numbering of the formulas,
it is straightforward to formalize most aspects of this statement.
There are primitive recursive functions ``Neg,'' ``Conj,'' and
``Imp'' encoding negation, conjunction, and implication (i.e.,
if \itc{q} is the G\"odel number of~\itc{Q} then Neg~\itc{q}
is the G\"odel number of the negation of~\itc{Q,} and so on),
and a primitive recursive relation~\itc{R} that relates \itc{i}
to~\itc{j} if and only if \itc{i} is the G\"odel number of a
proof of the sentence whose G\"odel number is~\itc{j.}
The only tricky part is the self-reference,
and this is achieved using the
usual (primitive recursive) ``diagonalization'' operator~\itc{D}:
\itc{D(m,n)} is the
G\"odel number of the sentence obtained by replacing the
free variable in the formula having G\"odel number \itc{m}
by the name of the number~\itc{n.}

Now let \itc{q$_1$} and~\itc{q$_2$} be the G\"odel numbers
of \Do\ and \Dt\ respectively, let $\not\equiv$ denote
exclusive or, and let \itc{P}[\itc{x}] abbreviate
$\exists$\itc{y}:\itc{yRx.}  \itc{P} stands for provable.
Then we can formulate the following formula with the free
variable~\itc{x}:

\bigskip
\centerline{(\Do\ \& $\sim$\itc{P}[\itc{D(x,x)} Imp \itc{q$_1$}])
$\not\equiv$ (\Dt\ \& $\sim$\itc{P}[(\itc{D(x,x)} Conj Neg \itc{q$_1$})
Imp \itc{q$_2$}]).}

Let \itc{h} be the G\"odel number of this formula,
and let \itc{S} be the sentence obtained by substituting
\itc{h} for~\itc{x,} i.e.,

\bigskip
\centerline{\hfill (\Do\ \& $\sim$\itc{P}[\itc{D(h,h)} Imp \itc{q$_1$}])
$\not\equiv$ (\Dt\ \& $\sim$\itc{P}[(\itc{D(h,h)} Conj Neg \itc{q$_1$})
Imp \itc{q$_2$}]).\hfill $(\dag)$}

\noindent
Then, by definition of~\itc{D,}
\itc{D(h,h)} is the G\"odel number of~\itc{S.}
The clincher is that \itc{D(h,h)} also appears on the right-hand side
of~$(\dag)$ exactly where we want it to appear.
Using ``\#'' to denote ``the name of the G\"odel number of'' we can
rewrite \itc{S} as

\bigskip
\centerline{(\Do\ \& $\sim$\itc{P}[\#(\itc{S} $\Rightarrow$ \Do)])
$\not\equiv$ (\Dt\ \& $\sim$\itc{P}[\#((\itc{S} \& $\sim$\Do)
$\Rightarrow$ \Dt)]).}

\noindent
This completes the formalization of statement~(c).
We can now imitate the students' argument to
show that \itc{S} is logically false, i.e., that
$\sim$\itc{S} is a tautology.  Using the definition of~\itc{S,}
we can prove

\bigskip
\centerline{\hfill(\itc{S} \& $\sim$\Do) $\Rightarrow$ \Dt.\hfill(1)}

Let \itc{a} be the G\"odel number of~(1).
By the nature of the relation~\itc{R,}
the provability of (1) implies~\itc{P(a).}
But observe that $\sim$\itc{P(a)} appears in the second disjunct in the
definition of~\itc{S.}  It follows that

\bigskip
\centerline{\hfill\itc{S}~$\Rightarrow$~\Do.\hfill(2)}

\noindent
The rest of the argument is
now clear: if \itc{b} is the G\"odel number of~(2), then
the provability of (2) implies \itc{P(b),} but
$\sim$\itc{P(b)} appears in the first disjunct of~\itc{S.}
Therefore~$\sim$\itc{S.}

Thus, although self-reference is not
illegitimate in all circumstances, it is illegitimate here because
this particular self-referential statement is self-contradictory.
Fitch's proof has a satisfying air of definitiveness,
and seems to vindicate~Shaw.

However, various authors have raised objections to this analysis.
The most important is that the proof does not give
any explanation for
why the teacher's announcement appears to be vindicated after the fact.
It appears to pin the blame on the teacher's
announcement instead of on the students,
and surely this cannot be correct.

A related objection rests on the observation that
if the teacher had not announced the exam
to the class but had simply decided in secret to give a surprise exam,
then no paradox would have occurred.
Therefore the trouble cannot be attributed solely to
the propositional content of the teacher's announcement;
the act of announcing it to the students must play a crucial role.
The purely logical analysis seems to ignore this.

These objections have convinced many to reject entirely the ``purely
logical'' approach, and
to propose a different, ``epistemological'' approach.

Before moving on to a discussion of the epistemological school,
however, I want to point out that
the objections \itc{can} be met.  For
example, the first objection indicates a misunderstanding of the
purely logical approach.  The conclusion of the logical analysis is
\itc{not} that the teacher's announcement is self-contradictory and is
the source of the paradox.  Rather, the conclusion is that \itc{in
order for the students to carry out their argument} that the teacher's
announcement cannot be fulfilled, they must \itc{interpret} the
teacher's announcement as saying something like~(c).
If the teacher intended (c) when making the
announcement, then it would be contradictory, and would
remain so after the examination.  However, a more reasonable assumption
is that the teacher's announcement, whatever it means, does not mean~(c),
and that therefore the students misinterpret the announcement when
they make their argument.  The announcement appears to be vindicated
afterwards, but the statement that is actually vindicated is
something like
``the students will be \itc{psychologically} surprised by
the exam,'' and such a statement does not permit the students' argument
to be carried out.  Similar observations are made in [4] and~[9].

As for the objection about the role of the act of making the announcement,
observe that the same sequence of words can have different meanings
depending on context, and that in the case of the teacher's announcement,
the public utterance of the sentence changes its propositional
content from ``there will be a surprise exam'' to something like ``there
will be a surprise exam in spite of the fact that I am now telling you
that there will be a surprise exam.''  The logical analysis therefore
\itc{does} take into account the act of making the announcement, albeit
implicitly, in its definition of the word ``surprise.''  Ignoring the
act of making the announcement would leave us stuck at~(a).

\sec The epistemological school

The purely logical approach is attractive to a mathematician
both because it shows exactly what problems arise from trying
to convert the paradoxical argument into a mathematical proof
and because it has connections to nontrivial theorems of logic.
However, it has one serious disadvantage: certain
aspects of the paradox---the act of announcing the exam, the
belief or disbelief that the students have in the announcement,
their assumption that they will remember the announcement during
the course of the week, and so on---are taken into account only
implicitly and not explicitly.  It is therefore natural to ask
if we can formalize the paradox in a way that lays bare these
``epistemic'' aspects.

Various epistemological formalizations
have been proposed in the literature;
we give just one here (taken from~[29])
to illustrate the idea.
As before, reduce the number of days to two
for simplicity; let ``1'' denote ``the exam occurs on the first day''
and let ``2'' denote ``the exam occurs on the second day.'' Let ``Ka''
denote ``on the eve of the first day the students will know'' and let
``Kb'' denote ``on the eve of the second day the students will know.''
The announcement can then be written

\bigskip

\def\nt{\mathord{\sim}}
\centerline{\hfill [1 $\Rightarrow$ $\nt$Ka 1] \&
  [2 $\Rightarrow$ ($\nt$Kb 2 \& Kb $\nt$1)] \& [1 $\vee$ 2].\hfill ($\ddag$)}

\noindent
We now introduce certain assumptions about knowledge and add them to
our list of rules of inference in our logic.

\item{KD:} If one knows \itc{A} \& \itc{B,} then one knows \itc{A}
and one knows \itc{B.}  Similarly, if one knows that \itc{A} implies
\itc{B} and one knows \itc{A,} then one knows \itc{B.}

\item{KI:} All logical truths are known.

\item{KE:} It is not possible to know something that is false.

We begin the argument with a lemma: Kb($\ddag$) $\Rightarrow$ $\nt$2;
remember that ``$\Rightarrow$'' here
encompasses our new rules of logic KD, KI, and KE.
Assume that Kb($\ddag$) is true.  By KD, it follows that
Kb[1 $\vee$ 2].  Now assume towards a contradiction that~2
is true, i.e., the exam is held on the last day.  From Kb($\ddag$)
and KE, ($\ddag$) follows, and 2 together with second conjunct of~($\ddag$)
implies $\nt$Kb 2 \& Kb $\nt$1; in particular, $\nt$Kb 2.
On the other hand, using KD, we deduce from Kb[1 $\vee$ 2]
and Kb $\nt$1 that Kb 2, a contradiction.
Thus, Kb($\ddag$) implies $\nt$2, and by KI we can infer
Ka[Kb($\ddag$) $\Rightarrow$ $\nt$2].

Now we can proceed with the crux of the argument:
deducing a contradiction from the assumption KaKb($\ddag$).
Assume KaKb($\ddag$).  From KD and Ka[Kb($\ddag$) $\Rightarrow$ $\nt$2]
it follows that Ka $\nt$2.
It is one of our logical truths (KE) that Kb($\ddag$) $\Rightarrow$ ($\ddag$),
so from KI we conclude that Ka[Kb($\ddag$) $\Rightarrow$ ($\ddag$)].
By KD and our assumption that KaKb($\ddag$), this implies
Ka($\ddag$) and in particular (by KD again) that Ka(1 $\vee$ 2).
Since we know that Ka $\nt$2, it follows from KD
that Ka~1.  But since Ka($\ddag$) is true,
($\ddag$) is true (by KE), and in particular its
first disjunct 1 $\Rightarrow$ $\nt$Ka 1 is true.  Then from Ka~1
we deduce 1 (from KE) and hence $\nt$Ka~1, a contradiction.

This shows that certain plausible assumptions
about knowledge---KI, KD, and KE, together with the assumption that
the students know that they will know the content of the announcement
throughout the week---are inconsistent.
Pointing out to the students that
they are making these internally inconsistent assumptions about
knowledge is enough to dissolve the paradox;
we do not necessarily have to decide which assumption is the ``wrong'' one.

It is still interesting, however,
to see if one of the assumptions appears to be
a particularly promising candidate for rejection.
Perhaps the most popular candidate
has been the assumption that
after hearing the announcement, the students ``know'' the content of the
announcement.  Those who maintain that we can never ``know'' things
by authority or that we can never ``know'' things about the future (at
least not with the same certainty that we can know many other things)
naturally find this approach attractive.
However, even those who are less skeptical have reason
to reject the assumption, because the statement that the students are
supposed to ``know'' is a statement that says something about the
students' \itc{inability} to ``know'' certain things.
For comparison, consider the statement,
``It is raining but John Doe does not know that it is raining.''
Clearly, John Doe cannot know the content
of this statement even if the statement is true and it is uttered in
his hearing by an extraordinarily reliable source.
This curious phenomenon is known as a ``Moore paradox''
or a ``blindspot,'' and the surprise exam paradox may be viewed as
simply a more intricate version of this situation.
The easiest way to see the connection is to reduce the
length of the week to one day, so that the announcement
becomes, ``There will be an exam tomorrow but you do not know that.''
This approach is essentially the one offered in
[3], [6], [18], [21], [22], [24], and~[29].

Others have argued that the assumption
KaKb($\ddag$) is plausible only
if one invokes the ``temporal retention principle''
(the students know that they will not forget the announcement
during the week) or ``Hintikka's KK principle'' (if one knows
something then one knows that one knows it), and that one or
both of these assumptions should be discarded.
I do not discuss this in detail here since I feel it
is of limited mathematical interest, but I mention
a brilliant variation of the paradox concocted by Sorensen~[28],
which suggests that rejecting these assumptions may be missing the point.

Exactly one of five students, Art, Bob, Carl, Don, and Eric, is to
be given an exam.  The teacher lines them up alphabetically so
that each student can see the backs of the students ahead
of him in alphabetical order but not the students after him.
The students are shown four silver stars and one gold star.
Then one star is secretly put on the back of each student.
The teacher
announces that the gold star is on the back of the student
who must take the exam, and that that student will be
surprised in the sense that he will not know he has been
designated until they break formation.  The students argue
that this is impossible; Eric cannot be designated because
if he were he would see four silver stars and would know that
he was designated.  The rest of the argument proceeds in
the familiar way.  The significance of this variation is
that in our preceding formalization we can let Ka mean ``Art knows''
and Kb mean ``Bob knows'' and then KaKb($\ddag$) appears to be
immediately plausible without reference to time or the KK
principle.  Thus, the problem remains even if those principles
are rejected;
see [28] and [16] for more discussion.

A very interesting variant of the epistemological approach,
that of Kaplan and Montague~[15], is a kind of
hybrid of the logical and epistemological schools.
They prove a theorem called ``the Paradox of the Knower''
that is reminiscent of
Tarski's theorem on the indefinability of the truth predicate.
Suppose we we have a first-order language and we
wish to introduce a knowledge predicate~\itc{K.}
There are certain reasonable-sounding conditions that we might
want to place on~\itc{K}:

\bigskip
\leftline{(A) \itc{K}(\#\itc{Q}) $\Rightarrow$ \itc{Q};}
\leftline{(B) (A) is known, i.e., \itc{K}(\#A);}
\leftline{(C) if \itc{Q} can be proved from \itc{P} and
\itc{K}(\#\itc{P}), then \itc{K}(\#\itc{Q}).}

\noindent
Unfortunately, these assumptions cannot be satisfied.
Using a diagonalization argument,
we can construct a sentence~\itc{S} such
that \itc{S}~$\equiv$~\itc{K}(\#($\sim$\itc{S})), and then
derive a contradiction by substituting this~\itc{S} for \itc{Q} and
A for \itc{P} in (A), (B), and~(C).
Thus, no such knowledge predicate is possible.

One might think at first that (C) is the dubious assumption since
certainly nobody knows all the logical consequences of what he knows,
but (C) can be weakened to the assumption that the logical
conclusion of a \itc{particular} explicitly given proof is known,
so the theorem is quite a strong one.
The Paradox of the Knower has inspired some sophisticated work in logic;
see [1], [2], or~[13].

\sec Game theory

Some authors have made the fascinating suggestion that
the surprise exam paradox
may be related to the iterated prisoner's dilemma.
The prisoner's dilemma is a two-player game in which
each player has the choice of either defecting or cooperating and must
make the choice without communicating with the other player and without
prior knowledge of the other player's choice.  If one
player defects and the other player cooperates, then the defector
enjoys a large payoff and the cooperator suffers a large loss.  If
both players defect then both payoffs are zero, and if both
players cooperate then they both earn a moderate payoff.
It is easy to show that each player has a dominant strategy
(i.e., one that is better than any other strategy
regardless of the opponent's strategy): to defect.

Intuitively, defection is the best choice because the prisoner's dilemma
is a ``one-shot'' game; there is no incentive for players to build up
a cooperative relationship since they are guaranteed never to meet
again.  This suggests considering the iterated prisoner's dilemma,
in which there are \itc{n}~rounds instead of just one (and the fact
that there are exactly \itc{n}~rounds is public knowledge).  The
payoffs in each round are as in the usual prisoner's dilemma, and
the two players are still not allowed to communicate with each other,
but at each round they do know and remember the results of all
previous rounds.  One might
think that in this case, occasional cooperation would be superior to
invariable defection---the idea being that a cooperative move in an
early stage, even if the opponent defects, encourages future
cooperation that counterbalances earlier losses.

Consider the following ``surprise examination''
argument that even in the \itc{n}-round
prisoner's dilemma, the optimal strategy is invariable defection.
The last round of an iterated prisoner's dilemma is identical to the
``one-shot'' prisoner's dilemma, since there is no hope of future
cooperation.  Hence the optimal last-round strategy is to defect.
But since defection in the last round is certain, there is no
incentive in the penultimate round to cooperate, for doing so
cannot possibly encourage future cooperation.  Thus, the optimal
strategy in the penultimate round is also defection.  Proceeding
by induction, we conclude that perfect players always defect.

The analogy between this argument and the standard surprise examination
argument is quite striking at first.
Indeed, Sorensen~[30] has argued that the two are really the same,
and has substantially revised his analysis of the surprise exam
as a result.
There is, however, an important disanalogy.
In the iterated prisoner's dilemma,
the conclusion about invariable defection is
\itc{counterintuitive,} but it does not lead to an explicit
\itc{contradiction.}
It is not difficult to adapt our argument to
give a fully rigorous mathematical proof that in the iterated prisoner's
dilemma, a Nash equilibrium is possible
only if both players defect in every round;
see [12, p.~166].
(A Nash equilibrium is a
situation in which if the strategies of all but one player are held fixed,
that player cannot do better by changing strategies.
One reason that the ``surprise exam''
argument we presented is not rigorous as it
stands is that the word ``optimal'' is imprecise.
A Nash equilibrium is a precise concept that
captures some---though not all---of the connotations
of the word ``optimal.'')
Therefore, I believe that the iterated prisoner's dilemma
is essentially distinct from the surprise examination paradox and
is not just a variant; see~[23].

Nevertheless, one might be able to exploit
the parallel between the surprise examination and the
iterated prisoner's dilemma to obtain some new ideas for game theory.
After all, cooperation is observed in the real world,
and this suggests that the usual mathematical model of the
iterated prisoner's dilemma might ignore some crucial point.
For some interesting ideas in this direction, see~[17].

Finally, I want to mention an unpublished idea of Karl Narveson
that illustrates how the surprise exam paradox can inspire new mathematics.
A teacher gives a quiz every week, with probability {\it p}$_1$
on Monday, {\it p}$_2$ on Tuesday, and so on.
The teacher's goal is to find a probability distribution that
maximizes the absolute value of the expected surprise when the quiz
is announced.  Here ``surprise'' is based on Shannon entropy,
so the surprise on Monday is log {\it p}$_1$,
the surprise on Tuesday is the log of the probability that the
exam occurs on Tuesday given that it has not occurred on Monday,
and so on until Friday, when the
quiz becomes a certainty and its announcement no longer comes as
a surprise.

Let {\it q}$_{n,m}$ be the
probability that the exam occurs on the \itc{n}th-to-the-last
day of an \itc{m}-day week given that it has not occurred on
any previous days, where \itc{n} ranges from zero to {\it m} $-$ 1.
As one can easily show, the optimal value of {\it q}$_{n,m}$
is independent of~\itc{m,} so we drop the second subscript.

Now set {\it s}$_0$~=~0.
Narveson has shown that {\it q}$_n$ is given by the mutual recursions
$$
\eqalign{\hbox{{\it q}$_n$} &= \hbox{exp({\it s}$_{n-1} - 1$)} \cr
\hbox{{\it s}$_n$} &= \hbox{{\it s}$_{n-1}$ $-$ {\it q}$_n$} \cr}$$
The {\it p}'s may then be recovered from the {\it q}'s.
For a five-day week, the probabilities for each of the five days are
about 0.1620, 0.1654, 0.1713, 0.1844, and 0.3169.

\sec Further reading

The literature contains a wide variety of other approaches to the
surprise examination paradox.
Cargile~[5] is the first paper in the literature to mention game theory.
Clark~[7] remarks that strictly mathematical analyses of the surprise exam
are rare in the literature and he tries to fill this gap.
Smullyan~[27] weaves G\"odel's theorem,
brainteasers, the surprise exam,
and other ``epistemic'' and ``doxastic'' paradoxes
into a delightful tapestry.
Some have seen connections between the surprise exam and the
sorites paradox (``removing one grain of sand from a heap of sand
leaves it a heap so zero grains of sand is still a heap'');
see [8], [26], and~[31].
A connection with the paradox of Schr\"odinger's cat is discussed in
[14] and~[19].

\sec Acknowledgments

The USENET newsgroup {\tt rec.puzzles} has played an important role in
my understanding of the surprise examination paradox.  I am indebted to
an unknown contributor to the {\tt rec.puzzles} archive who compiled
a list of about thirty references; reading these papers some five years
ago constituted my first serious attempt to understand the paradox.
I am grateful for the encouragement of Chris Cole, the maintainer
of the {\tt rec.puzzles} archive, and of Bill Taylor, whose enlightening
and witty correspondence has had a lot of influence on this paper.
Roy Sorensen, Dale Jacquette, and Martin Gardner have helped
with my literature search.
Thanks finally to Jim Propp for suggesting that I submit this paper
to the M{\sc ONTHLY}, and to all my friends who have patiently
endured my ramblings on the paradox over the years.

\sec References

{
\advance\baselineskip by -2pt
\font\rm=pplr at 8pt
\font\it=pplri at 8pt
\font\bf=pplb at 8pt
\parskip=2.5pt plus.5pt minus.5pt
\rm

\item{1.}
C. A. Anderson, The paradox of the knower, \itc{J. Phil.}\
{\bf 80} (1983), 338--355.

\item{2.}
N. Asher and H. Kamp, The knower's paradox and representational theories
of attitudes, in \itc{Theoretical Aspects of Reasoning About Knowledge:
Proceedings of the 1986 Conference, March 19--22, 1986, Monterey,
California,} ed.~J.~Y.~Halpern, Morgan Kaufmann, Los Altos, California
(1986), pp.~131--147.

\item{3.}
R. Binkley, The surprise examination in modal logic, \itc{J.
Phil.}\ {\bf 65} (1968), 127--136.

\item{4.}
J. Bosch, The examination paradox and formal prediction, \itc{Logique
et Analyse} {\bf 15} (1972), 505--525.

\item{5.}
J. Cargile, The surprise test paradox, \itc{J. Phil.}\
{\bf 64} (1967), 550--563.

\item{6.}
J. M. Chapman and R. J. Butler, On Quine's `so-called paradox,' \itc{Mind}
{\bf 74} (1965), 424--425.

\item{7.}
D. Clark, How expected is the unexpected hanging?
\itc{Math.\ Mag.}\ {\bf 67} (1994), 55--58.

\item{8.}
P. Dietl, The surprise examination, \itc{Educational Theory} {\bf 23}
(1973), 153--158.

\item{9.}
M. Edman, The prediction paradox, \itc{Theoria} {\bf 40} (1974), 166--175.

\item{10.}
H. B. Enderton, \itc{A Mathematical Introduction to Logic,}
Academic Press, New York (1972).

\item{11.}
F. Fitch, A Goedelized formulation of the prediction paradox,
\itc{Amer.\ Phil.\ Quart.}\ {\bf 1} (1964), 161--164.

\item{12.}
D. Fudenberg and J. Tirole, \itc{Game Theory,} MIT Press, Cambridge, MA
(1991).

\item{13.}
P. Grim, Operators in the paradox of the knower, \itc{Synthese} {\bf 94}
(1993), 409--428.

\item{14.}
J. M. Holtzman, A note on Schr\"odinger's cat and the unexpected
hanging paradox, \itc{British J. Phil.\ Sci.}\
{\bf 39} (1988), 397--401.

\item{15.}
D. Kaplan and R. Montague, A paradox regained, \itc{Notre Dame J.
Formal Logic} {\bf 1} (1960), 79--90.

\item{16.}
R. L. Kirkham, On paradoxes and a surprise exam, \itc{Philosophia}
{\bf 21} (1991), 31--51.

\item{17.}
R. C. Koons, Doxastic paradox and reputation effects in iterated games,
\itc{Theoretical Aspects of Reasoning About Knowledge:
Proceedings of the Fourth Conference (TARK 1992), March 22--25, 1992,
Monterey, California,} ed.~Y.~Moses, Morgan Kaufmann, San Mateo, California
(1992), pp.~60--72.

\item{18.}
I. Kvart, The paradox of surprise examination, \itc{Logique et Analyse}
{\bf 21} (1978), 337--344.

\item{19.}
J. G. Loeser, Three perspectives on Schr\"odinger's cat,
\itc{Amer.\ J.\ Physics} {\bf 52} (1984), 1089--1093;
letters and replies, {\bf 53} (1985), 937 and {\bf 54} (1986), 296--297.

\item{20.}
A. Margalit and M. Bar-Hillel, Expecting the unexpected, \itc{Philosophia}
{\bf 13} (1983), 263--288.

\item{21.}
T. H. O'Beirne, Can the unexpected \itc{never} happen? \itc{New Scientist}
{\bf 10} (1961), 464--465; letters and replies, 597--598.

\item{22.}
D. Olin, The prediction paradox resolved, \itc{Phil.\ Studies}
{\bf 44} (1983), 225--233.

\item{23.}
D. Olin, Predictions, intentions, and the prisoner's dilemma,
\itc{Phil.\ Quart.}\ {\bf 38} (1988), 111--116.

\item{24.}
W. V. O. Quine, On a so-called paradox, \itc{Mind} {\bf 62} (1953), 65--67.

\item{25.}
R. Shaw, The paradox of the unexpected examination, \itc{Mind}
{\bf 67} (1958), 382--384.

\item{26.}
J. W. Smith, The surprise examination on the paradox of the heap,
\itc{Phil.\ Papers} {\bf 13} (1984), 43--56.

\item{27.}
R. Smullyan, \itc{Forever Undecided: A Puzzle Guide to
G\"odel,} Knopf, New York (1987), parts I--V, particularly chapter~2.

\item{28.}
R. A. Sorensen, Recalcitrant versions of the prediction paradox,
\itc{Australasian J. Phil.}\ {\bf 69} (1982), 355--362.

\item{29.}
R. A. Sorensen, Conditional blindspots and the knowledge squeeze:
a solution to the prediction paradox, \itc{Australasian J.
Phil.}\ {\bf 62} (1984), 126--135.

\item{30.}
R. A. Sorensen, \itc{Blindspots,} Clarendon Press, Oxford (1988).

\item{31.}
T. Williamson, Inexact knowledge, \itc{Mind} {\bf 101} (1992), 217--242.

\item{32.}
P. Windt, The liar in the prediction paradox, \itc{Amer.\ Phil.\
Quart.}\ {\bf 10} (1973), 65--68.

\bigskip\bigskip
\parskip=0pt
\parindent=0pt
\obeylines
\it
Tellabs Research Center
One Kendall Square
Cambridge, MA 02139
tchowATalumDOTmitDOTedu
}

\vfill\eject
\rm

\footline={\hss\rm\folio\hss}
\parskip=10pt plus3pt minus3pt
\advance\parindent by 15pt
\def\itc#1{{\it #1\/}}
\hfuzz=3pt
\hyphenation{Philo-sophy}

\centerline{\bf The Surprise Examination or Unexpected Hanging Paradox:}
\centerline{\bf An Exhaustive Bibliography}
\medskip
\centerline{Timothy Y. Chow}
\centerline{\tt tchowATalumDOTmitDOTedu}
\bigskip

When writing~[Cho], I tried very hard to
compile an exhaustive bibliography for the surprise exam paradox.
I believe I came very close to succeeding,
but the M{\sc ONTHLY} refused to publish it,
saying that extensive bibliographies
are not appropriate for M{\sc ONTHLY} articles.
I am therefore making the bibliography freely available on the web instead.
I will try to update this bibliography continually;
if you know of any references that are not listed here,
please let me know!

The list below was last updated June~25, 2011.
(The last previous ArXiv update was made on September~6, 2005.)
Minor updates have been made to some entries:
[Clk], [Cra], [Frn], [Jac], [K-M], [Kea], [Loe], [Sob].
As one might expect, Google Scholar proved to be an invaluable
tool for discovering new references, including some such as
[Lan], [Ge1], [Pri], [Kaw], and [Chn]
that I had previously overlooked despite publication dates of
1986, 1999, 2000, 2002, and 2002 respectively.
At the same time, Google Scholar generated many hits
that I was not sure I wanted to include.
Some people had written substantial articles
or given talks on the paradox but had not (yet) formally published them;
in the end, I decided to exclude unpublished items from this bibliography
(though I did allow for dissertations such as [Mar]
and electronic publications such as [So13]).
Other papers were primarily concerned with some other topic
and made only passing reference to the paradox,
or used the paradox just as an illustrative example of something.
These papers I evaluated on a case-by-case basis
and made a personal judgment call;
I decided to include [Ald], [Cuo], [Da1], [Da2],
[Fas], [Gin], [Gra], [H-H], [Kl1], [Kl2], [L-K], and [WHR],
but omitted a number of others.
Specifically, I decided that the Kaplan--Montague paradox of the knower was
distinct from the surprise examination paradox,
and I even went back and deleted some entries
that I had included in previous versions of my bibliography
that were concerned only with the paradox of the knower.

Following my previous practice,
I included books that contained a substantial discussion
of the paradox even if most of the book was devoted to other topics;
in this category were [Cv2], [Cv3], [Kim], [So12], and~[T-H].
And of course,
I documented a steady stream of new publications on the paradox:
[Cv1], [DK1], [DK2], [F-B],
[Ge2], [G-C], [Hu], [Hua], [Kna], [K-R], [Lvy],
[\O GS], [Prr], [Sav], [S-W], [V\'az], [Wm1], [Wm2], and [Zha].
Of these, [Prr] (see also the accompanying commentary in~[\O GS])
is one of the most interesting because it is a previously
unpublished manuscript that probably dates from the early 1950's.

Finally, in addition to the acknowledgments that I made
in previous versions of this bibliography,
I would like to thank Sandy Lemberg and Arthur Benjamin
for notifying me about some of the new references above,
and Hans van Ditmarsch, Paul Franceschi, Lasse Burri Gram-Hansen,
Ran Raz, and John Williams for keeping me informed about their
own research on the paradox.

{\bf References}
\medskip

\advance\baselineskip by -2pt
\font\rm=pplr at 8pt
\font\it=pplri at 8pt
\font\bf=pplb at 8pt
\font\tt=cmtt8
\parskip=2.5pt plus.5pt minus.5pt
\rm

\item{[Ald]}
D. Aldous, Stopping times and tightness II, \itc{The Annals of
Probability} {\bf 17} (1989), 586--595.

\item{[Ale]}
P. Alexander, Pragmatic paradoxes, \itc{Mind} {\bf 59} (1950), 536--538.

\item{[Aud]}
R. Audi (ed.), \itc{The Cambridge Dictionary of Philosophy,}
Cambridge University Press, New York (1995), entries on
\itc{paradox} and \itc{unexpected examination paradox.}

\item{[Au1]}
A. K. Austin, On the unexpected examination, \itc{Mind} {\bf 78} (1969),
137.

\item{[Au2]}
A. K. Austin, The unexpected examination, \itc{Analysis} {\bf 39}
(1979), 63--64.

\item{[Aye]}
A. J. Ayer, On a supposed antinomy, \itc{Mind} {\bf 82} (1973), 125--126.

\item{[B-C]}
J. Bennett and J. Cargile, Reviews, \itc{Journal of Symbolic Logic}
{\bf 30} (1965), 101--103.

\item{[Ber]}
J. L. Berm\'udez, Rationality and the backwards induction argument,
\itc{Analysis} {\bf 59} (1999), 243--248.

\item{[Bin]}
R. Binkley, The surprise examination in modal logic, \itc{Journal of
Philosophy} {\bf 65} (1968), 127--136.

\item{[Blk]}
S. Blackburn, \itc{The Oxford Dictionary of Philosophy,} 
Oxford University Press, Oxford (1994), entry on \itc{prediction paradox.}

\item{[Bla]}
U. Blau, Vom Henker, vom L\"ugner und vom ihrem Ende,
\itc{Erkenntnis} {\bf 19} (1983), 27--44.

\item{[BB1]}
J. Blau and U. Blau, Epistemische Paradoxien I, \itc{Dialectica}
{\bf 49} (1995), 169--193.

\item{[BB2]}
J. Blau and U. Blau, Epistemische Paradoxien II, \itc{Dialectica}
{\bf 50} (1996), 167--182.

\item{[BBM]}
D. Borwein, J. M. Borwein and P. Mar\'echal, Surprise maximization,
\itc{American Mathematical Monthly} {\bf 107} (2000), 517--527.

\item{[Bos]}
J. Bosch, The examination paradox and formal prediction, \itc{Logique
et Analyse} {\bf 15} (1972), 505--525.

\item{[Bov]}
L. Bovens, The backward induction argument for the finite iterated
prisoner's dilemma and the surprise exam paradox, \itc{Analysis}
{\bf 57} (1997), 179--186.

\item{[Bun]}
B. Bunch, \itc{Mathematical Fallacies and Paradoxes,} Van Nostrand,
New York (1982), 34--37.

\item{[Car]}
J. Cargile, The surprise test paradox, \itc{Journal of Philosophy}
{\bf 64} (1967), 550--563.

\item{[Cas]}
J. Case, Paradoxes involving conflicts of interest, \itc{American
Mathematical Monthly} {\bf 107} (2000), 33--43.

\item{[Cv1]}
P. Cave, Reeling and A-reasoning: Surprise examinations and
Newcomb's tale, \itc{Philosophy} {\bf 79} (2004), 609--616.

\item{[Cv2]}
P. Cave, \itc{Can a Robot Be Human?  33 Perplexing Philosophy Puzzles,}
Oneworld (2009), Chapter~11.

\item{[Cv3]}
P. Cave, \itc{This Sentence is False: An Introduction to Philosophical
Paradoxes,} Continuum (2009), Chapter~2.

\item{[Cha]}
T. S. Champlin, Quine's judge, \itc{Philosophical Studies} {\bf 29}
(1976), 349--352.

\item{[C-B]}
J. M. Chapman and R. J. Butler, On Quine's `so-called paradox,' \itc{Mind}
{\bf 74} (1965), 424--425.

\item{[Chi]}
C. Chihara, Olin, Quine and the surprise examination, \itc{Philosophical
Studies} {\bf 47} (1985), 191--199.

\item{[Cho]}
T. Y. Chow, The surprise examination or unexpected hanging paradox,
\itc{American Mathematical Monthly} {\bf 105} (1998), 41--51.

\item{[Cla]}
D. Clark, How expected is the unexpected hanging?
\itc{Mathematics Magazine} {\bf 67} (1994), 55--58.

\item{[Clk]}
M. Clark, \itc{Paradoxes from A to Z,} Routledge, London, 2002,
entries on \itc{The Unexpected Examination} and
\itc{The Designated Student.}

\item{[Cl1]}
R. Clark, Pragmatic paradox and rationality, \itc{Canadian Journal of
Philosophy} {\bf 24} (1994), 229--242.

\item{[Cl2]}
R. Clark, Reflection and truth, in 
\itc{Theoretical Aspects of Reasoning About Knowledge:
Proceedings of the Fourth Conference (TARK 1992), March 22--25, 1992,
Monterey, California,} ed.~Y.~Moses, Morgan Kaufmann, San Mateo, California
(1992), 73--84.

\item{[Chn]}
D. H. Cohen, Informal logic \& the surprise exam,
\itc{Informal Logic} {\bf 22} (2002), TS 22--27.

\item{[Coh]}
L. J. Cohen, Mr.~O'Connor's ``pragmatic paradoxes,'' \itc{Mind} {\bf 59}
(1950), 85--87.

\item{[Cra]}
E. Craig (ed.), \itc{Routledge Encyclopedia of Philosophy,}
Routledge, New York (1998), entry on \itc{paradoxes, epistemic};
also in \itc{Concise Routledge Encyclopedia of Philosophy,}
Routledge, New York (2000), and
\itc{The Shorter Routledge Encyclopedia of Philosophy,} New York (2005).

\item{[Cuo]}
M. A. Cuonzo, How to solve paradoxes: A taxonomy and analysis
of solution-types, \itc{Cogency} {\bf 1} (2009), 9--21.

\item{[D-S]}
J. Dancy and E. Sosa (eds.), \itc{A Companion to Epistemology,}
Basil Blackwell Ltd., Cambridge, MA (1992),
entry on \itc{surprise examination paradox.}

\item{[Da1]}
E. Davis, A first-order theory of communicating first-order formulas,
\itc{Proceedings of the Ninth International Conference on Principles
of Knowledge Representation and Reasoning, June 2--5, 2004, Whistler,
Canada}, ed.~D.~Dubois, C.~A. Welty, and M.-A. Williams,
AAAI Press (2004), 235--245.

\item{[Da2]}
E. Davis, Knowledge and communication: A first-order theory,
\itc{Artificial Intelligence} {\bf 166} (2005), 81--139.

\item{[Die]}
P. Dietl, The surprise examination, \itc{Educational Theory} {\bf 23}
(1973), 153--158.

\item{[DK1]}
H. van Ditmarsch and B. Kooi, The secret of my success,
\itc{Synthese} {\bf 151} (2006), 201--232;
erratum (republication): {\bf 153} (2006), 339.

\item{[DK2]}
H. van Ditmarsch and B. Kooi, Een analyse van de Hangman Paradox
in dynamische epistemische logica, {\itc Algemeen Nederlands
Tijdschrift voor Wijsbegeerte} {\bf 97} (2005), 16--30.

\item{[Edm]}
M. Edman, The prediction paradox, \itc{Theoria} {\bf 40} (1974), 166--175.

\item{[Edw]}
P. Edwards (ed.), \itc{The Encyclopedia of Philosophy,}
Crowell Collier and MacMillan, Inc., New York (1967),
entry on \itc{logical paradoxes.}

\item{[E-F]}
G. W. Erickson and J. A. Fossa (eds.),
\itc{Dictionary of Paradox,}
University Press of America, Inc., Lanham, MD (1998),
entries on \itc{Hollis's paradox} and \itc{prediction paradox, the.}

\item{[Fal]}
N. Falletta, \itc{The Paradoxicon,}
Doubleday and Co., Inc., Garden City, NY (1983), 161--166.

\item{[Fas]}
M. Fasli, Towards a first-order approach for social agents:
Preliminary report, \itc{Proceedings of the Twelfth International
Florida Artificial Intelligence Research Society (FLAIRS) Conference,}
ed.~A.~Kumar and I.~Russell, 1999.

\item{[Fer]}
K. G. Ferguson, Equivocation in the surprise exam paradox, \itc{Southern
Journal of Philosophy} {\bf 29} (1991), 291--302.

\item{[F-B]}
J. L. Ferreira and J. Z. Bonilla, The surprise exam paradox,
rationality, and pragmatics: A simple game-theoretic analysis,
\itc{Journal of Economic Methodology} {\bf 15} (2008), 285--299.

\item{[Fit]}
F. Fitch, A Goedelized formulation of the prediction paradox,
\itc{American Philosophical Quarterly} {\bf 1} (1964), 161--164.

\item{[Frn]}
P. Franceschi,
Une analyse dichotomique pour le paradoxe de l'examen-surprise,
\itc{Philosophiques} {\bf 32} (2005), 399--421.
[Formerly titled: A straightforward reduction of the surprise
examination paradox to the sorites paradox.]

\item{[Fra]}
J. T. Fraser, Note relating to a paradox of the temporal order, in \itc{The
Voices of Time,} ed.~J.~T.~Fraser, University of Massachusetts Press,
Amherst (1981), 524--526.

\item{[Ful]}
J. S. Fulda, The paradox of the surprise test, \itc{The Mathematical
Gazette} {\bf 75} (1991), 419--421.

\item{[Gal]}
P. Galle, Another note on the ``surprise test'' puzzle, \itc{Informal Logic
Newsletter} {\bf 3} (1981), 21--22.

\item{[GaS]}
G. Gamow and M. Stern, The date of the hanging, in \itc{Puzzle-Math,}
Viking Press, New York (1958), 23--27.

\item{[Ga1]}
M. Gardner, A new prediction paradox, \itc{British Journal for the
Philosophy of Science} {\bf 13} (1962), 51.

\item{[Ga2]}
M. Gardner, A new paradox, and variations on it, about a man condemned to be
hanged, \itc{Scientific American} {\bf 208} (March 1963), 144--154;
reprinted in:
\itc{The Unexpected Hanging and Other Mathematical Diversions:
With a New Afterword and Expanded Bibliography,} University of Chicago
Press, Chicago (1991), chapter~1 and afterword;
\itc{The Colossal Book of Mathematics: Classic Puzzles,
Paradoxes, and Problems,}
W. W. Norton \& Company, New York (2001), chapter~43;
\itc{Martin Gardner's Mathematical Games}
[CD-ROM], Mathematical Association of America, 2005.

\item{[Ga3]}
M. Gardner, \itc{Further Mathematical Diversions,} Penguin Books (1977),
chapter I.

\item{[Ga4]}
M. Gardner, The erasing of Philbert the fudger, in \itc{Science Fiction
Puzzle Tales,} Clarkson Potter, New York (1981), puzzle~20.

\item{[Ga5]}
M. Gardner, \itc{aha! Gotcha: Paradoxes to Puzzle and Delight,}
W. H. Freeman and Company, San Francisco (1982), 26--27.

\item{[Ga6]}
M. Gardner, Again, how's that again? in \itc{Riddles of the Sphinx and
Other Mathematical Puzzle Tales,} Mathematical Association of America,
Washington, D.C. (1987), chapter~27.

\item{[Gea]}
J. Geanakoplos,
The hangman's paradox and Newcomb's
paradox as psychological games,
Cowles Foundation Discussion Paper \#1128,
Cowles Foundation, Yale University, 1996.

\item{[Ge1]}
J. Gerbrandy, Bisimulations on Planet Kripke,
Ph.D. thesis, University of Amsterdam, 1999.

\item{[Ge2]}
J. Gerbrandy, The surprise examination in dynamic epistemic logic,
\itc{Synthese} {\bf 155} (2007), 21--33.

\item{[GiS]}
I. Gilboa and D. Schmeidler,
Information dependent games: Can common sense be common knowledge?
\itc{Economics Letters} {\bf 27} (1988), 215--221.

\item{[Gin]}
H. Gintis, Rationality and its discontents
(review of \itc{Rational Decisions} by Ken Binmore),
\itc{The Economic Journal} {\bf 120} (2010), F162--F180.

\item{[Go1]}
L. Goldstein, Inescapable surprises and acquirable intentions,
\itc{Analysis} {\bf 53} (1993), 93--99.

\item{[Go2]}
L. Goldstein, Examining boxing and toxin,
\itc{Analysis} {\bf 63} (2003), 242--244.

\item{[G-C]}
L. Goldstein and P. Cave, A unified Pyrrhonian resolution of
the toxin problem, the surprise examination, and Newcomb's puzzle,
\itc{American Philosophical Quarterly} {\bf 45} (2008), 365--376.

\item{[Gra]}
R. Grafstein, Taking Dworkin to Hart: A positivist conception of
institutional rules, \itc{Political Theory} {\bf 11} (1983), 244-265.

\item{[Gui]}
S. Guiasu, Prediction paradox revisited, \itc{Logique et Analyse}
{\bf 30} (1987), 147--154.

\item{[Hal]}
N. Hall, How to set a surprise exam, \itc{Mind} {\bf 108} (1999), 647--703.

\item{[H-M]}
J. Y. Halpern and Y. Moses, Taken by surprise: the paradox of the
surprise test revisited, \itc{Journal of Philosophical Logic} {\bf 15}
(1986), 281--304.

\item{[Hrd]}
R. Hardin, \itc{Collective Action,} Johns Hopkins University Press,
Baltimore (1982), 145--150.

\item{[Har]}
C. Harrison, The unanticipated examination in view of Kripke's semantics
for modal logic, in \itc{Philosophical Logic,} ed.~J.~W.~Davis,
D.~J.~Hockney and K.~Wilson, D.~Reidel, Dordrecht (1969), 74--88.

\item{[H-H]}
J. Heather and D. Hill, I'm not signing that!
\itc{Proceedings of the First International Workshop on Formal Aspects
in Security and Trust (FAST 2003), September 8--9, 2003, Pisa,}
ed.~T.~Dimitrakos and F.~Martinelli (2003), 71--82.

\item{[Hel]}
N. S. K. Hellerstein, \itc{Delta, A Paradox Logic,} World Scientific,
Singapore (1997), 230--232.

\item{[Hil]}
M. Hild, Auto-epistemology and updating, \itc{Philosophical Studies}
{\bf 92} (1998), 321--361.

\item{[Ho1]}
M. Hollis, A paradoxical train of thought, \itc{Analysis} {\bf 44}
(1984), 205--206.

\item{[Ho2]}
M. Hollis, More paradoxical epistemics, \itc{Analysis} {\bf 46}
(1986), 217--218.

\item{[Hon]}
T. Honderich (ed.), \itc{The Oxford Companion to Philosophy,}
2nd ed., Oxford University Press, New York (2005),
entries on \itc{examination paradox} and \itc{prediction paradox.}

\item{[Hz1]}
J. M. Holtzman, An undecidable aspect of the unexpected hanging problem,
\itc{Philosophia} {\bf 17} (1987), 195--198.

\item{[Hz2]}
J. M. Holtzman, A note on Schr\"odinger's cat and the unexpected
hanging paradox, \itc{British Journal for the Philosophy of Science}
{\bf 39} (1988), 397--401.

\item{[Hu]}
Y.-L. Hu, Dui ``Turan yanxi beilun'' de fenxi (Analysis of the ``surprise
drill paradox''), \itc{Journal of Chongqing Institute of Technology}
(2007), no.~12, 81--85.

\item{[Hua]}
B. Huang, Egg, tiger, and hangman---Analyses on epistemic paradoxes,
\itc{Journal of Southwest University (Social Science Edition)}
(2010), no.~5.

\item{[HSW]}
C. Hudson, A. Solomon, R. Walker, A further examination of the `surprise
examination paradox,' \itc{Eureka} {\bf 32} (1969), 23--24.

\item{[H-B]}
P. Hughes and G. Brecht, \itc{Vicious Circles and Infinity: A Panoply
of Paradoxes,} Doubleday, Garden City, NY (1975), 41--61.

\item{[Jac]}
F. Jackson, The easy examination paradox, in
\itc{Analytical Philosophy in Comparative Perspective: Exploratory
Essays in Current Theories and Classical Indian Theories of Meaning
and Reference,} Springer (1984), 151--159;
reprinted with minor revisions in
F. Jackson, \itc{Conditionals,} Basil Blackwell, Oxford (1987), chapter~7.

\item{[Ja1]}
D. Jacquette, Review of Roy Sorensen's \itc{Blindspots, Journal of
Speculative Philosophy} {\bf 3} (1989), 218--223.

\item{[Ja2]}
D. Jacquette, A deflationary resolution of the surprise event paradox,
\itc{Iyyun: The Jerusalem Philosophical Quarterly} {\bf 41} (1992),
335--349.

\item{[Ja3]}
D. Jacquette, On the designated student and related induction paradoxes,
\itc{Canadian Journal of Philosophy} {\bf 24} (1994), 583--592.

\item{[Jan]}
C. Janaway, Knowing about surprises: a supposed antinomy revisited,
\itc{Mind} {\bf 98} (1989), 391--409.

\item{[JK1]}
B. Jongeling and T. Koetsier, A reappraisal of the hangman paradox,
\itc{Philosophia} {\bf 22} (1993), 299--311.

\item{[JK2]}
B. Jongeling and T. Koetsier, Blindspots, self-reference,
and the prediction paradox,
\itc{Philosophia} {\bf 29} (2002), 377--391.

\item{[Kal]}
M. G. Kalin, A tighter noose for the hangman, \itc{International Logic
Review} {\bf 22} (1980), 105--106.

\item{[Kan]}
S. Kanger, The paradox of the unexpected hangman, in \itc{Wright and Wrong:
Mini Essays in Honor of Georg Henrik von Wright on His Sixtieth Birthday,
June~14, 1976,} ed.~K.~Segerberg, Group in Logic and
Methodology of Real Finland (1976), 19--23.

\item{[K-M]}
D. Kaplan and R. Montague, A paradox regained, \itc{Notre Dame Journal
of Formal Logic} {\bf 1} (1960), 79--90.
Reprinted in \itc{Formal Philosophy: Selected Papers of Richard Montague,}
ed.~R.~H. Thomason, Yale University Press (1974), 271--285.

\item{[Kaw]}
Y. Kawano, G\"odel's incompleteness and forcing,
\itc{Algorithms in Algebraic Systems and Computation Theory,}
RIMS Kokyuroku {\bf 1268} (2002), 126--137.

\item{[Kea]}
J. Kearns, An illocutionary logical explanation of the surprise execution,
\itc{History and Philosophy of Logic} {\bf 20} (1999), 195--213.

\item{[K-E]}
J. Kiefer and J. Ellison, The prediction paradox again, \itc{Mind}
{\bf 74} (1965), 426--427.

\item{[Kim]}
S. O. Kimbrough, \itc{Agents, Games, and Evolution: Strategies at
Work and Play,} Taylor \& Francis (in press), Chapter~19.

\item{[Kin]}
M. Kinghan, A paradox derailed: reply to Hollis, \itc{Analysis}
{\bf 46} (1986), 20--24.

\item{[Ki1]}
R. L. Kirkham, The two paradoxes of the unexpected examination,
\itc{Philosophical Studies} {\bf 49} (1986), 19--26.

\item{[Ki2]}
R. L. Kirkham, On paradoxes and a surprise exam, \itc{Philosophia}
{\bf 21} (1991), 31--51.

\item{[Kl1]}
R. Klein, Under ``pragmatic'' paradoxes,
\itc{Yale French Studies} {\bf 66}, The Anxiety of Anticipation (1984),
91--109.

\item{[Kl2]}
R. Klein, The future of nuclear criticism,
\itc{Yale French Studies} {\bf 77}, Reading the Archive: On Texts and
Institutions (1990), 76--100; reprinted in
\itc{Yale French Studies} {\bf 97}, 50 Years of Yale French Studies:
A Commemorative Anthology, Part~2: 1980--1998 (2000), 78--102.

\item{[Kna]}
M. J. Knauff, The surprise examination paradox: A rejection of Quine,
and alternate solutions, \itc{The Yale Philosophy Review} {\bf 4}
(2008), 8--16.

\item{[Kra]}
M. H. Kramer, Another look at the problem of the unexpected examination,
\itc{Dialogue} {\bf 38} (1999), 491--501.

\item{[K-R]}
S. Kritchman and R. Raz, The surprise examination paradox and
the second incompleteness theorem, \itc{Notices of the American Mathematical
Society} {\bf 57} (2010), 1454--1458;
letter in response by S.~Alexander, {\bf 58} (2011), 527.

\item{[Kva]}
I. Kvart, The paradox of surprise examination, \itc{Logique et Analyse}
{\bf 21} (1978), 337--344.

\item{[Lan]}
F. Landman, Paradoxes of elimination, in \itc{Towards a Theory of
Information: The Status of Partial Objects in Semantics,}
GRASS 6, Foris (1986), 183--228.

\item{[Lee]}
B. D. Lee, The paradox of the surprise examination revisited,
in \itc{Circularity, Definition and Truth}, ed.~A.~Chapuis,
Indian Council of Philosophical Research (2000), 247--254.

\item{[Lei]}
J. Leiber, \itc{Paradoxes,}
Gerald Duckworth \& Co., Ltd., London (1993), 79--80.

\item{[Len]}
W. Lenzen, Die Paradoxie der \"uberraschenden \"Ubung: logische,
epistemologische und pragmatische Aspekte, \itc{Logique et Analyse}
{\bf 19} (1976), 267--284.

\item{[Lev]}
D. S. Levi, Surprise! \itc{Southern Journal of Philosophy}
{\bf 38} (2000), 447--464.

\item{[Lvy]}
K. Levy, The solution to the surprise exam paradox,
\itc{The Southern Journal of Philosophy} {\bf 47} (2009), 131--158.

\item{[L-K]}
J. Li and G. Kendall, Finite iterated prisoner's dilemma revisited:
Belief change and end-game effect, \itc{Proceedings of the Behavioral
and Quantitative Game Theory (BQGT '10): Conference on Future Directions,
Mary 14--16, 2010, Newport Beach, CA, USA} (2010), Article~48.

\item{[Lin]}
B. Linsky, Factives, blindspots and some paradoxes, \itc{Analysis}
{\bf 46} (1986), 10--15.

\item{[Loe]}
J. G. Loeser, Three perspectives on Schr\"odinger's cat,
\itc{American Journal of Physics} {\bf 52} (1984), 1089--1093;
letters and replies, P.~J. Bussey, {\bf 53} (1985), 937
and C.~Kacser, {\bf 54} (1986), 296--297.

\item{[Lyo]}
A. Lyon, The prediction paradox, \itc{Mind} {\bf 68} (1959), 510--517.

\item{[Mac]}
J. J. Macintosh,
Aquinas and Ockham on time, predestination and the unexpected examination,
\itc{Franciscan Studies} {\bf 55} (1998), 181--220.

\item{[Mar]}
A. Marcoci, The surprise examination paradox in dynamic epistemic logic,
M.Sc. thesis, University of Amsterdam, 2010.

\item{[M-B]}
A. Margalit and M. Bar-Hillel, Expecting the unexpected, \itc{Philosophia}
{\bf 13} (1983), 263--288.

\item{[McC]}
J. McClelland, Epistemic logic and the paradox of the surprise examination,
\itc{International Logic Review} {\bf 3} (1971), 69--85.

\item{[M-C]}
J. McClelland and C. Chihara, The surprise examination paradox,
\itc{Journal of Philosophical Logic} {\bf 4} (1975), 71--89.

\item{[Med]}
B. Medlin, The unexpected examination, \itc{American Philosophical
Quarterly} {\bf 1} (1964), 66--72; corrigenda,~333.

\item{[Mel]}
B. Meltzer, The third possibility, \itc{Mind} {\bf 73} (1964), 430--433.

\item{[M-G]}
B. Meltzer and I. J. Good, Two forms of the prediction paradox, \itc{British
Journal for the Philosophy of Science} {\bf 16} (1965), 50--51.

\item{[Met]}
U. Metschl, Eine kleine \"Uberraschung f\"ur Gehirne im Tank:
Eine skeptische Notiz zu einem antiskeptischen Argument,
\itc{Zeitschrift f\"ur philosophische Forschung} {\bf 43} (1989), 519--527.

\item{[M-H]}
J.-J. Ch.\ Meyer and W. van der Hoek,
A modal constrastive logic: The logic of `but,'
\itc{Annals of Mathematics and Artificial Intelligence}
{\bf 17} (1996), 291--313.

\item{[Nai]}
M. A. Nait Abdallah, Logic programming of some mathematical paradoxes,
in \itc{Fundamentals of Computation Theory: Proceedings, International
Conference FCT '89, Szeged, Hungary, August 21--25, 1989,}
ed.~J.~Csirik, J.~Demetrovics, and F.~G\'ecseg,
Lecture Notes in Computer Science, vol.~380,
Springer-Verlag, New York (1989), 349--361.

\item{[Ner]}
G. Nerlich, Unexpected examinations and unprovable statements, \itc{Mind}
{\bf 70} (1961), 503--513.

\item{[Nie]}
H. A. Nielsen, A note on the ``surprise test'' puzzle, \itc{Informal Logic
Newsletter} {\bf 2} (1979), 6--7.

\item{[OB1]}
T. H. O'Beirne, Can the unexpected \itc{never} happen? \itc{New Scientist}
{\bf 10} (1961), 464--465; letters and replies, 597--598.

\item{[OB2]}
T. H. O'Beirne,
\itc{Puzzles and Paradoxes: Fascinating Excursions in Recreational
Mathematics,} Dover Press, New York (1984), chapter~11.

\item{[OCa]}
M. J. O'Carroll, Improper self-reference in classical logic and the
prediction paradox, \itc{Logique et Analyse} {\bf 10} (1967), 167--172.

\item{[OCo]}
D. J. O'Connor, Pragmatic paradoxes, \itc{Mind} {\bf 57} (1948), 358--359.

\item{[\O GS]}
P. \O hrstr\o m, L.~.B. Gram-Hansen, and U. Sandborg-Petersen,
Time and knowledge: Some reflections on Prior's analysis of
the paradox of the prisoner, \itc{Synthese,} to appear.
DOI: 10.1007/s11229-011-0048-y

\item{[Ol1]}
D. Olin, The prediction paradox resolved, \itc{Philosophical Studies}
{\bf 44} (1983), 225--233.

\item{[Ol2]}
D. Olin, The prediction paradox: resolving recalcitrant variations,
\itc{Australasian Journal of Philosophy} {\bf 64} (1986), 181--189.

\item{[Ol3]}
D. Olin, On a paradoxical train of thought, \itc{Analysis} {\bf 46}
(1986), 18--20.

\item{[Ol4]}
D. Olin, On an epistemic paradox, \itc{Analysis} {\bf 47} (1987), 216--217.

\item{[Ol5]}
D. Olin, Predictions, intentions, and the prisoner's dilemma,
\itc{Philosophical Quarterly} {\bf 38} (1988), 111--116.

\item{[Ol6]}
D. Olin, \itc{Paradox,} McGill-Queen's University Press,
Montreal \& Kingston (2003), chapter~3.

\item{[Par]}
G. Pareti, I paradossi pragmatici: bibliografia, \itc{Rivista di Filosofia}
{\bf 69} (1978), 170--174.

\item{[Pit]}
V. Pittioni, Das Vorhersageparadoxon, \itc{Conceptus} {\bf 17} (1983),
88--92.

\item{[Pop]}
K. R. Popper, A comment on the new prediction paradox, \itc{British
Journal for the Philosophy of Science} {\bf 13} (1962), 51.

\item{[Pos]}
C. J. Posy, Epistemology, ontology, and the continuum,
in \itc{The Growth of Mathematical Knowledge,} ed.~E. Grosholz,
Synthese Library, vol.~289, Kluwer, Dordrecht (2000), 199--219.

\item{[Pou]}
W. Poundstone, \itc{Labyrinths of Reason: Paradox, Puzzles, and the
Frailty of Knowledge,} Anchor Press/Double\-day, New York (1988), chapter~6.

\item{[Pri]}
G. Priest, The logic of backwards inductions,
\itc{Economics and Philosophy} {\bf 16} (2000), 267--285.

\item{[Prr]}
A. N. Prior, The paradox of the prisoner in logical form,
\itc{Synthese,} to appear.  DOI: 10.1007/s11229-011-9947-z

\item{[Qui]}
W. V. O. Quine, On a so-called paradox, \itc{Mind} {\bf 62} (1953), 65--67;
reprinted in \itc{The Ways of Paradox, and Other Essays,} Harvard
University Press, Cambridge, MA (1976), chapter~2.

\item{[Reg]}
D. H. Regan, \itc{Utilitarianism and Cooperation,} Clarendon Press, Oxford
(1980), page~74.

\item{[Res]}
N. Rescher, \itc{Paradoxes: Their Roots, Range, and Resolution,}
Carus Publishing Co., Chicago and La Salle (2001), 112--114.

\item{[Sai]}
R. M. Sainsbury, \itc{Paradoxes,} 2nd edition,
Cambridge University Press, Cambridge (1995), sections 4.2--4.4.

\item{[Sav]}
M. vos Savant,
Ask Marilyn, \itc{PARADE,} December 5, 2004.

\item{[Sck]}
F. Schick, Surprise, self-knowledge, and commonality,
\itc{The Journal of Philosophy} {\bf 97} (2000), 440--453.

\item{[Sch]}
J. Schoenberg, A note on the logical fallacy in the paradox of the
unexpected examination, \itc{Mind} {\bf 75} (1966), 125--127.

\item{[S-W]}
B. Schumacher and M. Westmoreland,
Reverend Bayes takes the unexpected examination,
\itc{Math Horizons} (September 2008), 26--27.

\item{[Scr]}
M. Scriven, Paradoxical announcements, \itc{Mind} {\bf 60} (1951), 403--407.

\item{[Sha]}
S. C. Shapiro, A procedural solution to the unexpected hanging and sorites
paradoxes, \itc{Mind} {\bf 107} (1998), 751--761.

\item{[Shp]}
R. A. Sharpe, The unexpected examination, \itc{Mind} {\bf 74} (1965), 255.

\item{[Shw]}
R. Shaw, The paradox of the unexpected examination, \itc{Mind}
{\bf 67} (1958), 382--384.

\item{[Sla]}
B. H. Slater, The examiner examined, \itc{Analysis} {\bf 35} (1974), 49--50.

\item{[Smi]}
J. W. Smith, The surprise examination on the paradox of the heap,
\itc{Philosophical Papers} {\bf 13} (1984), 43--56.

\item{[Smu]}
R. Smullyan, \itc{Forever Undecided: A Puzzle Guide to
G\"odel,} Knopf, New York (1987), parts I--V, particularly chapter~2.

\item{[Sob]}
E. Sober, To give a surprise exam, use game theory,
\itc{Synthese} {\bf 115} (1998), 355--373;
addendum available on the author's website at
{\tt http://philosophy.wisc.edu/sober/surp-add.pdf}

\item{[So1]}
R. A. Sorensen, Recalcitrant versions of the prediction paradox,
\itc{Australasian Journal of Philosophy} {\bf 69} (1982), 355--362.

\item{[So2]}
R. A. Sorensen, Conditional blindspots and the knowledge squeeze:
a solution to the prediction paradox, \itc{Australasian Journal of
Philosophy} {\bf 62} (1984), 126--135.

\item{[So3]}
R. A. Sorensen, The bottle imp and the prediction paradox,
\itc{Philosophia} {\bf 15} (1986), 421--424.

\item{[So4]}
R. A. Sorensen, A strengthened prediction paradox, \itc{Philosophical
Quarterly} {\bf 36} (1986), 504--513.

\item{[So5]}
R. A. Sorensen, Blindspotting and choice variations of the prediction
paradox, \itc{American Philosophical Quarterly} {\bf 23} (1986),
337--352.

\item{[So6]}
R. A. Sorensen, The bottle imp and the prediction paradox,~II,
\itc{Philosophia} {\bf 17} (1987), 351--354.

\item{[So7]}
R. A. Sorensen, \itc{Blindspots,} Clarendon Press, Oxford (1988).

\item{[So8]}
R. A. Sorensen, The vagueness of knowledge, \itc{Canadian Journal of
Philosophy} {\bf 17} (1987), 767--804.

\item{[So9]}
R. A. Sorensen, The earliest unexpected class inspection, \itc{Analysis}
{\bf 53} (1993), 252.

\item{[So10]}
R. A. Sorensen, Infinite ``backward'' induction arguments, \itc{Pacific
Philosophical Quarterly} {\bf 80} (1999), 278--283.

\item{[So11]}
R. A. Sorensen, Formal problems about knowledge,
in \itc{The Oxford Handbook of Epistemology,}
ed.~P.~K. Moser, Oxford University Press, New York (2002), 539--568.

\item{[So12]}
R. A. Sorensen, Paradoxes of rationality,
in \itc{The Oxford Handbook of Rationality,}
ed.~A.~R. Mele and P.~Rawling, Oxford University Press (2004), 267--268.

\item{[So13]}
R. A. Sorensen, Epistemic paradoxes,
in \itc{The Stanford Encyclopedia of Philosophy} (online), Spring 2009 edition,
{\tt http://plato.stanford.edu}

\item{[Sta]}
M. Stack, The surprise examination paradox, \itc{Dialogue} {\bf 16}
(1977), 207--212.

\item{[Ste]}
I. Stewart, Paradox lost, \itc{Manifold} {\bf 10} (1971), 19--20;
D. Woodall replies, 21.

\item{[Tha]}
M. Thalos, Conflict and co-ordination in the aftermath of oracular
statements, \itc{The Philosophical Quarterly} {\bf 47} (1997), 212--226.

\item{[T-H]}
T. Tymoczko and J. Henle, \itc{Sweet Reason: A Field Guide to
Modern Logic,} Springer (1999), 420--428.

\item{[V\'az]}
M. V\'azquez, Knowledge, information and surprise,
\itc{tripleC} {\bf 7} (2009), 194--201.

\item{[V-E]}
V. Vieru and V. Enache, The surprise-test paradox: A formal study,
\itc{Studies in Informatics and Control} {\bf 8} (1999), 317--327.

\item{[Wtb]}
R. Weintraub, Practical solutions to the surprise examination paradox,
\itc{Ratio (New Series)} {\bf 8} (1995), 161--169.

\item{[Wei]}
P. Weiss, The prediction paradox, \itc{Mind} {\bf 61} (1952), 265--269.

\item{[WHR]}
H. Wietze de Haan, W. H. Hesselink, and G. R. Renardel de Lavalette,
Knowledge-based asynchronous programming, \itc{Fundamenta Informaticae}
{\bf 63} (2004), 259--281.

\item{[Wm1]}
J. Williams, A simple solution to the surprise exam paradoxes,
\itc{Research Collection School of Social Sciences} (2004), Paper~26.

\item{[Wm2]}
J. N. Williams, The surprise exam paradox: Disentangling two
\itc{reductios}, \itc{Journal of Philosophical Research} {\bf 32} (2007),
67--94.

\item{[Wi1]}
T. Williamson, Review of Roy Sorensen's \itc{Blindspots, Mind} {\bf 99}
(1990), 137--140.

\item{[Wi2]}
T. Williamson, Inexact knowledge, \itc{Mind} {\bf 101} (1992), 217--242.

\item{[Wi3]}
T. Williamson, \itc{Knowledge and Its Limits,}
Oxford University Press, New York (2000), chapter~6.

\item{[Win]}
P. Windt, The liar in the prediction paradox, \itc{American Philosophical
Quarterly} {\bf 10} (1973), 65--68.

\item{[Woo]}
D. R. Woodall, The paradox of the surprise examination, \itc{Eureka}
{\bf 30} (1967), 31--32.

\item{[Wr1]}
M. J. Wreen, Surprising the examiner, \itc{Logique et Analyse} {\bf 26}
(1983), 177--190.

\item{[Wr2]}
M. J. Wreen, Passing the bottle, \itc{Philosophia} {\bf 15} (1986), 427--444.

\item{[W-S]}
C. Wright and A. Sudbury, The paradox of the unexpected examination,
\itc{Australasian Journal of Philosophy} {\bf 55} (1977), 41--58.

\item{[Wri]}
J. A. Wright, The surprise exam: prediction on last day uncertain,
\itc{Mind} {\bf 76} (1967), 115--117.

\item{[Wu]}
K. J. Wu, Believing and disbelieving, in \itc{The Logical Enterprise,}
ed.~A.~R. Anderson et al., Yale University Press, New Haven (1975),
211--219.

\item{[Zha]}
T.-S. Zhang, On hangman paradox, \itc{Academic Journal of Jinyang}
(2006), no.~1.

\bye